\documentclass[12pt,reqno]{amsart}
\usepackage{amsmath,amssymb,amscd,latexsym,amsthm,mathrsfs,comment, color}
\usepackage[unicode]{hyperref}
\usepackage{hypbmsec}
\newtheorem{Thm}{Theorem}
\newtheorem{Con}{Conjecture}
\newtheorem{Lem}{Lemma}

\newtheorem{cor}{Corollary}
\theoremstyle{remark}

\begin{document}

\title[The Erd\H{o}s--Moser equation revisited]%
{Moser's mathemagical work on the equation $1^k+2^k+\ldots+(m-1)^k=m^k$}
\author{Pieter Moree} 
\address{Max-Planck-Institut f\"ur Mathematik, Vivatsgasse 7, D-53111 Bonn, Germany}
\email{moree@mpim-bonn.mpg.de}
\date{\today}
\dedicatory{In memory of Alf van der Poorten {\rm (1942-2010)}}
\begin{abstract}
If the equation of the title has an integer solution with $k\ge 2$,
then $m>10^{10^6}$. Leo Moser showed this in 1953 by amazingly
elementary methods. With the hindsight of more than 50 years his
proof can be somewhat simplified. We give a further proof
showing that Moser's result can be derived from a von Staudt-Clausen type theorem. 
Based on more recent developments concerning this equation, we derive
a new result using the divisibility properties of numbers in the sequence
$\{2^{2e+1}+1\}_{e=0}^{\infty}$. In the final section we show that certain
Erd\H{o}s-Moser type equations arising in a recent paper of Kellner can be solved
completely.  
\end{abstract}

\maketitle

\section{Introduction}
\noindent In this paper we are interested in non-trivial solutions, that is 
solutions with $k\ge 2$, of the equation
\begin{equation}
\label{EME}
1^k+2^k+\ldots+(m-2)^k+(m-1)^k=m^k.
\end{equation}
The conjecture that such solutions do not exist was formulated
around 1950 by Paul Erd\H{o}s in a letter to Leo Moser.  For $k=1$, one has
the solution $1+2=3$ (and no further solutions). From now on
we will assume that $k\ge 2$.
Leo Moser \cite{Moser}
established the following theorem in 1953. 
\begin{Thm} 
\label{Leo}
{\rm (Leo Moser, 1953)}.
If
$(m,k)$ is a solution of \eqref{EME}, then
$m>10^{10^6}$.
\end{Thm}
His result has since been improved. 
Butske et al.~\cite{graphs} have shown by computing rather
than estimating certain quantities in Moser's original proof that $m>1.485\cdot 10^{9\,321\,155}$.
By proceeding along these lines this bound cannot be substantially improved.
 Butske et al.~\cite[p. 411]{graphs} expressed
the hope that new insights will eventually make it possible to reach the benchmark
$10^{10^7}$. 

The main purpose of this paper is to make Moser's remarkable proof of Theorem \ref{Leo} better known, and with the 
hindsight and technological developments of more than 50 years, to give an even cleaner version of Moser's proof. This is contained in Section \ref{sec:second}.\footnote{A large part of the material in
Section \ref{sec:second} is copied verbatim from Moser's paper.}
Moreover, we obtain the following refinement of Moser's result.
\begin{Thm}
\label{main}
Suppose that $(m,k)$ is a solution of {\rm (\ref{EME})} with $k\ge 2$, then\\
{\rm 1)} $m > 1.485\cdot 10^{9\,321\,155}$.\\  
{\rm 2)} $k$ is even, $m\equiv 3~({\rm mod~}8)$, $m\equiv \pm 1~({\rm mod~}3)$;\\
{\rm 3)} $m-1$, $(m+1)/2$, $2m-1$ and $2m+1$ are all square-free.\\
{\rm 4)} If $p$ divides at least one of the integers in (3), then $p-1|k$.\\
{\rm 5)} The number $(m^2-1)(4m^2-1)/12$ is square-free and has at least $4\,990\,906$ 
prime factors.
\end{Thm}
In fact, Moser proved (3) and (4) of Theorem \ref{main} and weaker versions of parts (2) and (5). Readers interested in the shortest (currently known) proof of Theorem \ref{main} are referred to Moree \cite{rabbits}.  The deepest result used to prove Theorem \ref{main} is Lemma \ref{Smod}. Using a binomial identity due to 
Pascal (1654) a reproof of Lemma \ref{Smod} was given recently by MacMillan and Sondow
\cite{MS}. To wit,
had Blaise Pascal's computing machine from 1642,  the 
Pascaline,\footnote{The Pascaline was originally developed for tax collecting purposes!} worked like a modern computer,  then Theorem \ref{main} could have been already proved in 1654.

In Section \ref{sec:comparison} we compare our alternative proof with Moser's original proof.

In Section \ref{sec:reproof} we give a more systematic proof of Moser's result, which uses
a variant of the von Staudt-Clausen theorem.\footnote{The proof given in
Section \ref{sec:reproof} is implicit in Moree's \cite{Oz} with $a=1$.} The relevance of this result for the study
of the Erd\H{o}s-Moser equation was first pointed out in 1996 by Moree \cite{Oz} who used 
the result to show that the Moser approach can also be used to study the equation
$1^k+2^k+\ldots+(m-1)^k=am^k$ and $a\ge 1$ an integer. An improvement of the main result of \cite{Oz} will be 
presented in Section \ref{eight}.

The reader might wonder which other techniques have been brought to bear for the study of (\ref{EME}).
Such techniques include Bernoulli numbers, considering the equation modulo prime powers, analysis (taking $k$
to be a real, rather than an integer) and continued fraction methods.  There is an extensive literature on the more general equation 
$$1^k+\ldots+(m-1)^k=y^n,~n\ge 2,$$ see, e.g., 
Bennett et al. \cite{BGP}. That work incorporates several further techniques. However, those results do not appear to have any implications for the study of (\ref{EME}).
In Section \ref{sec:MRU}, we give a taste of what can be done using Bernoulli numbers and considering
(\ref{EME}) modulo prime powers. The main result here is Theorem 1 of \cite{MRU}. We give a weakened
(far less technical) version of this, namely Lemma \ref{helpful}. Using that result and a heuristic
assumption on the behavior of $S_r(a)$, a heuristic argument validating the Erd\H{o}s-Moser
conjecture can be given (\cite[Section 6]{MRU-1}).

In Section \ref{rec:recent}, we consider implications for (\ref{EME}) based on analytic methods, and in 
particular the recent work of Gallot, Moree and Zudilin \cite{GMZ} who obtained the benchmark $10^{10^7}$ and further improved this to  
$10^{10^9}$ by computing $3\cdot 10^9$ digits of $\log 2$.

Section \ref{sec:new} is the most original part of the paper.  Results on divisors of numbers of the
form $2^{2e+1}+1$ are used to show that if $(m,k)$ is a solution of (\ref{EME}) such
that $m+2$ is only composed of primes $p$ satisfying $p\equiv 5,7~({\rm mod~}8)$, then
$m\ge 10^{10^{16}}$.

In the final two sections we consider the Erd\H{o}s-Moser variants
$$1^k+2^k+\ldots+(m-1)^k=am^k, {\rm ~respectively~}a(1^k+2^k+\ldots+(m-1)^k)=m^k$$ (with
$a\ge 1$ a fixed integer) and show that
the latter equation (arising in a recent paper of Kellner \cite{K}) can be solved completely for
infinitely many integers $a$.

This paper is in part scholarly and in part research.  Leo Moser (1921-1970) was a mathematician of the problem solver type.  For bibliographic information the reader
is referred to the MacTutor History of Mathematics archive \cite{Moser2} or Wyman \cite{Wyman}.

\section{Moser's proof revisited}
\label{sec:second}
Let $S_r(n)=\sum_{j=0}^{n-1}j^r$. In what follows we assume that
\begin{equation}
\label{M1}
S_k(m)=m^k,~k\ge 2,
\end{equation}
which corresponds to a non-trivial solution of (\ref{EME}). Throughout this note 
$p$ will be used to indicate primes.
\begin{Lem}
\label{Smod}
Let $p$ be a prime.
We have $$S_r(p)\equiv \epsilon_r(p)~({\rm mod~}p),$$
where
$$\epsilon_r(p)=
\begin{cases}
-1 & if~p-1|r;\\
 0 & otherwise.
\end{cases}
$$
\end{Lem}
\noindent {\it Proof}. Let $g$ be a primitive root modulo $p$. In case $p-1\nmid r$ we have
$$S_r(p)\equiv \sum_{j=0}^{p-2}(g^j)^r\equiv {g^{r(p-1)}-1\over g^r-1}~({\rm mod~}p),$$
and the numerator is divisible by $p$. In case $p-1|r$, we find by Fermat's Little
Theorem that $S_r(p)\equiv p-1\equiv -1~({\rm mod~}p)$ as desired. \qed\\

Another proof using only Lagrange's theorem on roots of polynomials over $\Bbb Z/p\Bbb Z$ 
can be given; see Moree \cite{rabbits}. The most elementary proof presently known is due to
MacMillan and Sondow \cite{MS} and is based on Pascal's identity (1654), valid for $n\ge 0$ and
$a\ge 2$:
$$\sum_{k=0}^n {n+1\choose k}S_k(a)=a^{n+1}-1.$$
A further proof can be given using the polynomial identity 
$$X^{p-1}-1\equiv \prod_{j=1}^{p-1}(X-j)~({\rm mod~}p)$$ and
Newton's identities expressing power sums in elementary symmetric polynomials.
\begin{Lem}
\label{pripower}
In case $p$ is an odd prime or in case $p=2$ and $r$ is even, we have
$S_r(p^{\lambda+1})\equiv pS_r(p^{\lambda})~({\rm mod~}p^{\lambda+1})$.
\end{Lem}
\noindent {\it Proof}. Every $0\le j<p^{\lambda+1}$ can be uniquely written as 
$j=\alpha p^{\lambda}+\beta$ with $0\le \alpha <p$ and $0\le \beta<p^{\lambda}$.
Hence we obtain by invoking the binomial theorem
$$S_r(p^{\lambda+1})=\sum_{\alpha=0}^{p-1}\sum_{\beta=0}^{p^{\lambda}-1}(\alpha p^{\lambda}+\beta)^r
\equiv
p\sum_{\beta=0}^{p^{\lambda}-1}\beta^r+rp^{\lambda}\sum_{\alpha=0}^{p-1}\alpha\sum_{\beta=0}^{p^{\lambda}-1}\beta^{r-1}~({\rm mod~}p^{2\lambda}) .$$
Since the first sum equals $S_r(p^{\lambda})$, and $2\sum_{\alpha=0}^{p-1}\alpha=p(p-1)\equiv 0~({\rm mod~}p)$,
the result follows. \qed\\

\noindent {\it Proof of Theorem} \ref{main}. Suppose that $p|m-1$, then using 
Lemma \ref{Smod} we infer that
\begin{equation}
\label{startje}
S_k(m)=\sum_{i=0}^{(m-1)/p-1}\sum_{j=1}^p(j+ip)^k\equiv {m-1\over p}S_k(p)\equiv {m-1\over p}\epsilon_k(p)~({\rm mod~}p).
\end{equation}
On the other hand $m\equiv 1~({\rm mod~}p)$, so that by (\ref{M1}) we must have
\begin{equation}
\label{M3}
{m-1\over p}\cdot \epsilon_k(p)\equiv 1~({\rm mod~}p).
\end{equation}
Hence $\epsilon_k(p)\not\equiv 0~({\rm mod~}p)$, so that from the definition of $\epsilon_k(p)$
it follows that $\epsilon_k(p)=-1$, and
\begin{equation}
\label{M4}
p|m-1 {\rm ~implies~}p-1|k.
\end{equation}
Thus (\ref{M3}) can be put in the form
\begin{equation}
\label{M5}
{m-1\over p}+1\equiv 0~({\rm mod~}p),
\end{equation}
or
\begin{equation}
\label{M6}
m-1\equiv -p~({\rm mod~}p^2).
\end{equation}
We claim that $m-1$ must have an odd prime divisor $p$, and that hence by (\ref{M4}),
$k$ must be even. 
It is easy to see that $m-1>2$.
If $m-1$ does not have an odd prime divisor, then $m-1=2^e$ for some $e\ge 2$. However,
by (\ref{M6}) we see that $m-1$ is square-free. This contradiction shows that
$m-1$ has indeed an odd prime factor $p$.

We now multiply together all congruences of the type (\ref{M5}), that is one for each prime
$p$ dividing $m-1$. Since $m-1$ is square-free, the resulting modulus is $m-1$. Furthermore,
products  containing  two or more distinct prime factors of the form $(m-1)/p$ will be
divisible by $m-1$. Thus we obtain
\begin{equation}
\label{M7}
(m-1)\sum_{p|m-1}{1\over p}+1\equiv 0~({\rm mod~}m-1),
\end{equation}
or 
\begin{equation}
\label{M8}
\sum_{p|m-1}{1\over p}+{1\over m-1}\equiv 0~({\rm mod~}1).
\end{equation}
We proceed to develop three more congruences, similar to (\ref{M8}), which when combined
with (\ref{M8}) lead to the proof of part 1. Equation (\ref{M1}) can be written in the
form
\begin{equation}
\label{M9}
S_k(m+2)=2m^k+(m+1)^k.
\end{equation}
Using Lemma \ref{Smod} and the fact that $k$ is even, we obtain as before
\begin{equation}
\label{M10}
p|m+1 {\rm ~implies~}p-1|k,
\end{equation}
and
\begin{equation}
\label{M11}
{m+1\over p}+2\equiv 0~({\rm mod~}p).
\end{equation}
{}From (\ref{M11}) it follows that no odd prime appears with exponent greater than
one in the prime factorization of $m+1$. The prime 2 (according to H. Zassenhaus `the oddest of primes'), requires  
special attention. If we inspect (\ref{EME}) with modulus 4 and use the fact
that $k$ is even, then we find that $m+1\equiv 1$ or $4~({\rm mod~}8)$. Now let us assume
that we are in the first case, and we let $2^f||m$ (that is $2^f|m$ and $2^{f+1}\nmid m$).
Note that $f\ge 3$. 
By an argument similar to that given in (\ref{startje}) we infer that 
$S_k(m+1)\equiv {m\over 2^f}S_k(2^f)~({\rm mod~}2^f)$.
Using Lemma
\ref{pripower}, we see that $S_k(m+1)\equiv {m\over 2^f}S_k(2^f)\equiv 2^{f-1}~({\rm mod~}2^f)$,
contradicting $S_k(m+1)=2m^k\equiv 0~({\rm mod~}2^f)$. Thus $m+1$ contains 2 exactly to the second power and hence
(\ref{M11}) can be put in the form
\begin{equation}
\label{M12}
{m+1\over 2p}+1\equiv 0~({\rm mod~}p).
\end{equation}
We multiply together all congruences of type (\ref{M12}). The modulus then becomes
$(m+1)/2$. Further, any term involving two or more distinct factors
${m+1\over 2p}$ will be divisible by ${m+1\over 2}$, so that on simplification we obtain
\begin{equation}
\label{M13}
\sum_{p|m+1}{1\over p}+{2\over m+1}\equiv 0~({\rm mod~}1).
\end{equation}
We proceed to find two similar equations to (\ref{M13}). Suppose that
$p|2m-1$, and let $t={1\over 2}({2m-1\over p}-1)$. Clearly $t$ is an integer,
and $m-1=tp+{p-1\over 2}$. We have $a^k=(-a)^k$ since $k$ is even so that
$2S_k({p+1\over 2})\equiv S_k(p)~({\rm mod~}p)$ and hence, by Lemma \ref{Smod},
$$S_k({p+1\over 2})\equiv {\epsilon_k(p)\over 2}~({\rm mod~}p).$$
It follows that
\begin{equation}
\label{M14}
S_k(m)\equiv \sum_{i=0}^{t-1}\sum_{j=1}^{p-1}(j+ip)^k+\sum_{i=1}^{(p-1)/2}i^k
\equiv (t+{1\over 2})\epsilon_k(p)~({\rm mod~}p).
\end{equation}

On the other hand $1\equiv (2m-1+1)^k\equiv (2m)^k~({\rm mod~}p)$, hence
$m^k\not\equiv 0~({\rm mod~}p)$, so that (\ref{M1}) and (\ref{M14}) imply
$\epsilon_k(p)\ne 0$. Hence $p-1|k$, and by Fermat's little theorem
$m^k\equiv 1~({\rm mod~}p)$. Thus (\ref{M1}) and (\ref{M14}) yield 
$-(t+{1\over 2})\equiv 1~({\rm mod~}p)$. Replacing $t$ by its value
and simplifying we obtain
\begin{equation}
\label{M15}
{2m-1\over p}+2\equiv 0~({\rm mod~}p).
\end{equation}
Since $2m-1$ is odd, (\ref{M15}) implies that $2m-1$ is square-free. 
Multiplying congruences of the type (\ref{M15}), one for each of the $r$ prime
divisors of $2m-1$, yields
$$2^{r-1}\Big((2m-1)\sum_{p|2m-1}{1\over p}+2\Big)\equiv 0~({\rm mod~}2m-1).$$
Since the modulus $2m-1$ is odd, this gives
\begin{equation}
\label{M16}
\sum_{p|2m-1}{1\over p}+{2\over 2m-1}\equiv 0~({\rm mod~}1).
\end{equation}
Finally we obtain a corresponding congruence for primes $p$ dividing $2m+1$, namely (\ref{M18}) 
below. For this
purpose we write (\ref{M1}) in the form
\begin{equation}
\label{M17}
S_k(m+1)=2m^k.
\end{equation}
Suppose $p|2m+1$. Set $v={1\over 2}({2m+1\over p}-1)$. Clearly $v$ is an integer. We
have $m=pv+{p-1\over 2}$ and find $S_k(m+1)\equiv (v+{1\over 2})\epsilon_k(p)~({\rm mod~}p)$.
From this and (\ref{M17}) it is easy to infer that $\epsilon_k(p)=-1$, and so $v+{1\over 2}\equiv -2~({\rm mod~}p)$.
We conclude that
$$p|2m+1 {\rm ~implies~}p-1|k.$$
Replacing $v$ by its value and simplifying, we obtain 
$${2m+1\over p}+4\equiv 0~({\rm mod~}p).$$
Note that this implies that $2m+1$ is square-free.
Reasoning as before we obtain
\begin{equation}
\label{M18}
\sum_{p|2m+1}{1\over p}+{4\over 2m+1}\equiv 0~({\rm mod~}1).
\end{equation}
If we now add the left hand sides of (\ref{M8}), (\ref{M13}), (\ref{M16}) and (\ref{M18}), we
get an integer, at least 4. 
By an argument similar to that showing $2\nmid m$, we show that $3\nmid m$ (but in this case we
use Lemma \ref{pripower} with $p=3$ and $3^{\lambda}||m$ and
the fact that $k$ must be even).
No prime $p>3$ can divide more than one of the integers
$m-1$, $m+1$, $2m-1$ and $2m+1$. Further, 
since $m\equiv 3~({\rm mod~}8)$ and $3\nmid m$,
2 and 3 divide precisely two of these integers.
We infer that $M_1=(m-1)(m+1)(2m-1)(2m+1)/12$ is a square-free integer. We deduce that
\begin{equation}
\label{M19}
\sum_{p|M_1}{1\over p}+{1\over m-1}+{2\over m+1}+{2\over 2m-1}+{4\over 2m+1}\ge 4-{1\over 2}-{1\over
3}=3{1\over 6}
\end{equation}
One checks that (\ref{M16}) has no solutions with $m\le 1\,000$. 
Thus (\ref{M19}) yields (with $\alpha=3.16$)
$\sum_{p|M_1}{1\over p}>\alpha$.
{}From this it follows that if 
\begin{equation}
\label{M20}
\sum_{p\le x}{1\over p}<\alpha, 
\end{equation}
then
$m^4/3>M_1>\prod_{p\le x}p$ and hence
\begin{equation}
\label{Pbound}
m > 3^{1/4}e^{\theta(x)/4},
\end{equation}
with $\theta(x)=\sum_{p\le x}\log p$, the Chebyshev $\theta$-function. Since
for example (\ref{M20}) is satisfied with $x=1\,000$, we find that $m>10^{103}$ and
infer from (\ref{M19}) that we can take $\alpha=3{1\over 6}-10^{-100}$ in (\ref{M20}).
Next one computes (using a computer algebra package, say PARI) the largest prime $p_k$ such that $\sum_{p\le p_k}{1\over p}<3{1\over 6}$, with
$p_1,p_2,\ldots$ the consecutive primes. Here one finds that $k=4\,990\,906$ and
$$\sum_{i=1}^{4\,990\,906}{1\over p_i}=3.166\,666\,658\,810\,172\,858\,4<3{1\over 6}-10^{-9}.$$ This completes 
the proof of part 1 of the theorem; the remaining parts of the theorem have been proven along the way.
\qed\\

\noindent {\tt Remark 1}. Since for a solution of (\ref{EME}), $(m^2-1)(4m^2-1)/12$ has at 
least 4\,990\,906 distinct prime factors, it is perhaps reasonable to expect that each of 
the factors $m-1$, $m+1$, $2m-1$
and $2m+1$ must have many distinct prime factors. Brenton and Vasiliu \cite{BV}, using
the bound given in part 1 of Theorem \ref{main},  showed 
that $m-1$ has
at least 26 prime factors. Gallot et al. \cite{GMZ} increased this, 
using Theorem \ref{YVES}, to 33.\\

\noindent {\tt Remark 2}. Moser considered (\ref{EME}) modulo $m-1,~m+1,2m-1$ and $2m+1$. Sondow
and MacMillan \cite{SM2} considered the equation also modulo $(m-1)^2$ and obtained some further
information (this involves the Fermat quotient).

\section{Comparison of the proof with Moser's}
\label{sec:comparison}
In this section we compare and contrast the proof of Theorem \ref{main} with Moser's proof of Theorem \ref{Leo}.  

Moser used only Lemma \ref{Smod}, not Lemma \ref{pripower}.  Consequently, he concluded that either $m\equiv 3~({\rm mod~}8)$ or $m\equiv 0~({\rm mod~}8)$. In the
first case we followed his proof but in the second case one has to note that we 
cannot use
(\ref{M13}). Letting $M_2=(m-1)(2m-1)(2m+1)$ we get from
(\ref{M8}), (\ref{M16}), (\ref{M18})
\begin{equation}
\label{M25}
\sum_{p|M_2}{1\over p}+{1\over m-1}+{2\over 2m-1}+{4\over 2m+1}>3-{1\over 3}
\end{equation}
However, since $2\nmid M_2$,
(\ref{M25}) is actually a stronger condition on $m$ than is
(\ref{M19}).\\
\indent The idea to use $3\nmid m$, leading to a slight improvement for the bound on
$m$, is taken from Butske et al. \cite{graphs} and not present in Moser's proof. (Actually they
consider the cases $3\nmid m$ and $3|m$ separately. We show that only $3\nmid m$ can occur.)\\
\indent By using some prime number estimates from Rosser, Moser deduces that
(\ref{M20}) holds with $x=10^7$ and $\alpha=3.16$. In his argument he claims that by direct
computation one sees that (\ref{M20}) holds with $x=1\,000$ and $\alpha=2.18$. This is not true (as pointed
out to me by Buciumas and Havarneanu).  %Reference?%  
However, replacing 2.18 by 2.2 in 
Moser's equation (21) one sees that his proof still remains valid. The present day possibilities
of computers allow us to proceed by direct computation, rather than 
to resort to prime number estimates as Moser was forced to do.\\
\indent The advantage of the proof given in Section \ref{sec:second} is that it shows,
in contrast to Moser's proof and Butske et al.'s variation thereof, that {\it every} 
non-trivial solution satisfies the crucial inequality (\ref{M19}).

\section{A second proof using a von Staudt-Clausen type theorem}
\label{sec:reproof}
In this section we show that Moser's four formulas (\ref{M8}), (\ref{M13}), (\ref{M16}) and (\ref{M18})
can be easily
derived from the following theorem. Indeed, using it a fifth formula can be  derived, namely (\ref{five})
below.
\begin{Thm} {\rm (Carlitz-von Staudt, 1961)}.
\label{thm:clausen}
Let $r,y$ be positive integers. Then
\begin{equation}
\label{staudt}
S_r(y)=\sum_{j=1}^{y-1}j^r=\begin{cases}
0~({\rm mod~}{y(y-1)\over 2}) &{\rm~if~}r{\rm ~is~odd};\cr
-\sum_{p-1|r,~p|y}{y\over p}~({\rm mod~}y) &{\rm ~otherwise}.
\end{cases}
\end{equation}
\end{Thm}
\noindent  Carlitz \cite{Carlitz} gave a proof of Theorem \ref{thm:clausen} using finite
differences and stated that the result is due
to von Staudt. In the case $r$ is odd, he claims that $S_r(y)/y$ is an integer, which
is not always true (it is true though that $2S_r(y)/y$ is always an integer). The 
author \cite{Canada} gave a proof of a generalization to sums of powers in arithmetic progression 
using the theory of primitive
roots. Kellner \cite{Kellner2} gave a reproof for even $r$ only) using Stirling numbers of 
the second kind. For the easiest proof known 
and some further applications of the Carlitz-von Staudt theorem, we refer the reader 
to Moree \cite{rabbits}.\\

\noindent {\it Second proof of Theorem} \ref{main}. We will apply Theorem \ref{thm:clausen} with $r=k$.\hfil\break
\indent In case $k$ is odd, we find by combining (\ref{staudt}) (with $y=m$) with (\ref{EME}) and using the
coprimality of $m$ and $m-1$ that $m=2$ or $m=3$, but these cases are easily excluded. Therefore
$k$ must be even.\\
\indent Take $y=m-1$. Then, using (\ref{EME}), the left hand side of (\ref{staudt}) simplifies to
$$S_k(m-1)=1^k+2^k+\ldots+(m-2)^k=m^k-(m-1)^k\equiv  1~({\rm mod~}m-1).$$
We get from (\ref{staudt}) that
\begin{equation}
\label{Pi1}
\sum_{p|m-1,~p-1|k}{(m-1)\over p}+1\equiv 0~({\rm mod~}m-1).
\end{equation}
Suppose there exists $p|m-1$ such that $p-1\nmid k$. Reducing both sides
modulo $p$, we get $1\equiv 0~({\rm mod~}p)$. This
contradiction shows that in (\ref{Pi1}) the condition $p-1|k$ can be dropped, and thus
we obtain (\ref{M7}). {}From (\ref{M7}) we see that $m-1$ must be square-free and also
we obtain (\ref{M8}).\\
\indent Take $y=m$. Then using (\ref{EME}) and $2|k$ we infer from (\ref{staudt}) that
\begin{equation}
\label{five}
\sum_{p-1|k,~p|m}{1\over p}\equiv 0~({\rm mod~}1).
\end{equation}
Since a sum of reciprocals of distinct primes can never be a positive integer,
we infer that the sum in (\ref{five}) equals zero and hence conclude that if $p-1|k$,
then $p\nmid m$. We conclude for example that $(6,m)=1$. Now on considering (\ref{EME}) with
modulus 4 we see that $m\equiv 3~({\rm mod~}8)$.\\
\indent Take $y=m+1$. Then using (\ref{EME}) and the fact that $k$ is even, the left hand side 
of (\ref{staudt}) simplifies to 
$$S_k(m+1)=S_k(m)+m^k=2m^k\equiv 2~({\rm mod~}m+1).$$
We obtain 
$$\sum_{p|m+1,~p-1|k}{(m+1)\over p}+2\equiv 0~({\rm mod~}m+1),$$
and by reasoning as in the
case $y=m-1$, it is seen that $p|m+1$ implies $p-1|k$, and thus (\ref{M13}) is obtained. {}From (\ref{M13}) and
$m\equiv 3~({\rm mod~}8)$, we derive that $(m+1)/2$ is square-free.\\
\indent Take $y=2m-1$. On noting that 
$$S_k(2m-1)=\sum_{j=1}^{m-1}(j^k+(2m-1-j)^k)\equiv 2S_k(m)\equiv 2m^k~({\rm mod~}2m-1),$$
we find that
\begin{equation}
\label{Pi2}
\sum_{p|2m-1,~p-1|k}{(2m-1)\over p}+2m^k\equiv 0~({\rm mod~}2m-1).
\end{equation}
Since $m$ and $2m-1$ are coprime, we infer that if $p|2m-1$, then $p-1|k$, $m^k\equiv 1~({\rm mod~}p)$
and furthermore that $2m-1$ is square-free. It follows from the
Chinese remainder theorem that $2m^k\equiv 2~({\rm mod~}2m-1)$, 
and hence from (\ref{Pi2}) we obtain (\ref{M16}).\\
\indent Take $y=2m+1$. Noting that
$$S_k(2m+1)=\sum_{j=1}^m(j^k+(2m+1-j)^k)\equiv 2S_k(m+1)=4m^k~({\rm mod~}2m+1)$$
and proceeding as in the case $y=2m-1$, we obtain (\ref{M18}) and the square-freeness of
$2m+1$.
To finish the proof we proceed as in Section \ref{sec:second} just below (\ref{M18}). \qed\\

\noindent With some of the magic behind the four Moser identities revealed, the reader might be well tempted to derive
further identities. A typical example would start from
\begin{equation}
4^k-1^k-2^k-3^k\equiv -\sum_{p-1|k\atop p|m-4}{(m-4)\over p}~({\rm mod~}m-4).
\end{equation}
For simplicity let us assume that $m\equiv 2({\rm mod~}3)$. We have $(6,m-4)=1$.
For this to lead to a further equation, we need the left hand side to be a constant
modulo $m-4$.
If we could infer that $p|m-4$ implies $p-1|k$, then the left hand side would equal $-2~({\rm mod~}m-4)$,
and we would be in business. (For the reader familiar with the Carmichael function $\lambda$, this can
be more compactly formulated as $\lambda(m-4)|k$.)
Unfortunately a problem is caused by the fact that the left hand side
could be divisible by $p$. Thus all we seem to obtain is that if $m\equiv 2~({\rm mod~}3)$, and
$\lambda(m-4)|k$ or $4^k-1^k-2^k-3^k$ and $m-4$ are coprime, then 
$$\sum_{p|m-4}{1\over p}-{2\over m-4}\equiv 0~({\rm mod~}1).$$

\noindent In Section \ref{sec:new}, we will see that if we replace $m-4$ by $m+2$ we can do a little better, the
reason being that in this case, $2^{k+1} +1$ appears on the left hand side, and numbers of this form have only a restricted set of possible
prime factors.

\section{Bernoulli numbers and a cascade process}
\label{sec:MRU}
Recall that the Bernoulli numbers $B_k$ are defined by the power series
$${t\over e^t-1}=\sum_{k=0}^{\infty}{B_kt^k\over k!}.$$ They are rational numbers and can be written
as $B_k=U_k/V_k$, with $(U_k,V_k)=1$. 
One has $B_0=1,~B_1=-1/2$ and $B_{2j+1}=0$ for $j\ge 1$.
By the von Staudt-Clausen theorem we can take for $k\ge 2$ even
$V_k=\prod_{p-1|k}p$. The Kummer congruences state that if $k$ and $r$ are even and
$k\equiv r\not\equiv 0~({\rm mod~}p-1)$, then $B_k/k\equiv B_r/r~({\rm mod~}p)$.
A prime $p$ will be called {\it regular} if it does not divide any
of the numerators $U_k$ with $k$ even and $\le p-3$. Otherwise it is said to be {\it irregular}. The first
few irregular primes are $37,~59,~67,~101,\ldots$.\\
\indent The power sum $S_r(n)$
can be expressed using Bernoulli numbers. One has, see e.g. \cite[(2.1)]{Murty},
$$S_r(n)=\sum_{j=0}^r {r\choose j}B_{r-j}{n^{j+1}\over j+1}.$$
Voronoi in 1889, see, e.g., \cite[Theorem 2.8]{Murty}), proved that if $k$ is even and $\ge 2$, then 
$V_rS_r(n)\equiv U_kn~({\rm mod~}n^2)$. {}From this
result we infer that for a solution $(k,m)$ of (\ref{EME}) we must have $m|U_k$ and thus in particular 
$\nu_p(U_k)\ge \nu_p(m)$, where we put $\nu_p(m)=f$ if $p^f||m$. By a more elaborate analysis Moree et
al. \cite{MRU} improved this to $\nu_p(B_k/k)\ge 2\nu_p(m)$. It shows (by the
von Staudt-Clausen theorem) that if $p|m$, then $p-1\nmid k$ (a 
conclusion we already reached using identity (\ref{five})). Invoking the Kummer congruences we 
then obtain the following
result.
\begin{Lem}
Let $(k,m)$ be a solution of {\rm (\ref{EME})} with $k\ge 2$ and even. If $p|m$, then $p$ is irregular.
\end{Lem}
\indent Let us call a pair $(r,p)$ with $p$ a regular prime and $2\le r\le p-3$ even, {\it helpful} if
for every $a=1,\ldots,p-1$ we have $S_r(a)\not\equiv a^r~({\rm mod~}p)$.
\begin{Lem}
\label{helpful}
If $(r,p)$ is a helpful pair, and $(k,m)$ a solution of {\rm (\ref{EME})} with $k$ even, then 
we have $k\not\equiv r~({\rm mod~}p-1)$.
\end{Lem}
{\it Proof}. Assume that $k\equiv r~({\rm mod~}p-1)$. By the previous lemma we must have $p\nmid m$. Now write $m=m_0p+b$. Thus $1\le b\le p-1$.
We have, modulo $p$, $S_k(m)\equiv S_r(m)\equiv m_0S_r(p)+S_r(b)\equiv S_r(b)$. Thus if (\ref{EME}) is
satisfied, we must have $S_r(b)\equiv b^r~({\rm mod~}p)$. By the definition of a helpful pair this
is impossible. \qed\\

\noindent Since $2|k$, and $(2,5)$ is a helpful pair, we infer that $4|k$. Since $(2,7)$ and $(4,7)$ are helpful pairs, it follows
that $6|k$. {}From $4|k$ and the fact that $(4,17)$ and $(12,17)$ are helpful pairs, it follows that $8|k$. We
thus infer that $24|k$. It turns out that this process can be continued to deduce that more and more small prime
factors must divide $k$; for
a detailed account with many tables see \cite{MRU-1}. Given an irregular prime $p$ and $2\le r\le p-3$ even, one would heuristically expect 
that it is helpful with probability $(1-1/p)^{p-1}$ which tends to $1/e$, assuming that the values
$S_r(a)$ are randomly distributed modulo $p$; this is supported by current numerical data.  \\
\indent Moree et al.~\cite{MRU}, using good pairs (of which the helpful pairs are a special case),
showed that
$N_1:=\operatorname{lcm}(1,2,\ldots,200)$ divides $k$. Kellner 
\cite{Kellner} showed in 2002 that also all primes $200<p<1000$ divide~$k$. Actually Moree et al. 
\cite[p.~814]{MRU}
proved
a slightly stronger result which combined with Kellner's shows that $N_2\mid k$ with
$$N_2=2^8\cdot 3^5\cdot 5^4\cdot 7^3\cdot 11^2\cdot 13^2\cdot 17^2\cdot 19^2\cdot 23\cdots
997>5.7462\cdot 10^{427}.$$
An heuristic argument can be given suggesting that if $L_v:={\rm lcm}(1,2,\ldots,v)$ divides $k$, with tremendously
high likelihood we can infer that $L_w$ divides $k$, where $w$ is the smallest prime power not dividing $L_v$. 
It suffices that $v\ge 11$.
To deduce that $k$ is divisible by 24 is delicate, but once one has $L_v|k$, there is an explosion of further 
helpful pairs one can use to establish divisibility of $k$ by a larger integer.
To add the first prime power $w$ not dividing $L_v$, one needs a number of helpful pairs that is roughly 
linear in $v$, whereas an exponential number (in $v$) is available.
However, the required computation time increases sharply with $w$. This heuristic argument is the most convincing known to the author in support of the Erd\H{o}s-Moser conjecture; details may be found in the 
extended version of \cite{MRU}, \cite{MRU-1}.\\
\indent Given a fixed integer $a$, one can try the same approach to study the equation $S_k(m)=am^k$. Again one sees that
once one manages to infer for example that $120|k$, one can show that there must be larger and larger divisors. For many $a$, however, this
`cascade process' does not seem to `take off' and it remains unknown whether all solutions with $k\ge 2$, if any, satisfy
$120|k$\, for example.\\
\indent If $n=\prod_i p_i^{e_i}$ denotes the canonical prime factorization of $n$, then $\Omega(n)=\sum e_i$ is
the total number of prime divisors of $n$. Urbanowicz \cite{U} proved a result which implies that given
an arbitrary $t$, there exists an integer $m_t$ such that if $(k,m)$ is a solution of (\ref{EME}) with $k\ge 2$
and $m\ge m_t$, then $\Omega(k)\ge t$.

\section{The analytical approach and continued fractions}
\label{rec:recent}
\indent Comparing $S_k(m)$ with the appropriate integrals, it is easy to see that the ratio $k/m$ must
be bounded. A more refined approach gives
$$S_k(m)={(m-1)^k\over 1-e^{-(k+1)/(m-1)}}(1+O({1\over \sqrt{m}})).$$
On equating the left-hand side to $m^k$ and using $(1-1/m)^m=\exp(-1+O(m^{-1}))$, one concludes that
as $m\rightarrow \infty$, we have
$${k\over m}=\log 2 + O({1\over \sqrt{m}}).$$
By a rather more delicate analysis Gallot et al. \cite{GMZ} obtain that for $m>10^9$ one has
$${k\over m}=\log 2\Big(1-{3\over 2m}-{C_m\over m^2}\Big),~{\rm~where~}0<C_m<0.004.$$
As a corollary this gives that if $(k,m)$ is a solution of (\ref{EME}) with $k\ge 2$ and even, then
$2k/(2m-3)$ is a convergent $p_j/q_j$ of $\log 2$ with $j$ even. 
This approach was first explored in 1976 by Best and
te Riele \cite{Best} in their attempt to solve a related conjecture of
Erd\H{o}s; see also Guy
\cite[D7]{Guy}.  The main result of \cite{GMZ} reads
as follows, where
given $N\ge 1$, we define
$$
{\mathcal{P}}(N)=\{p:p-1\mid N\}\cup\{p:\text{$3$ is a primitive root modulo $p$}\}.
$$
\begin{Thm}
\label{main3}
Let $N\ge 1$ be an arbitrary integer. Let
$$\frac{\log 2}{2N} = [a_0,a_1,a_2,\dots] = a_0 + \cfrac{1}{a_1 + \cfrac{1}{a_2 + \cdots}}$$
be the (regular) continued fraction of
$({\log 2})/(2N)$, with $p_i/q_i = [a_0,a_1,\dots,a_i]$ its $i$-th convergent.

Suppose that the integer pair $(m,k)$ with $k\ge 2$ satisfies \eqref{EME}
with $N\mid k$.
Let $j=j(N)$ be the smallest even integer such that:\\
{\rm a)} $a_{j+1}\ge 180N-2$;\\
{\rm b)} $(q_j,6)=1$;\\
{\rm c)} $\nu_p(q_j)=\nu_p(3^{p-1}-1)+\nu_p(N)+1$
for all primes $p\in {\mathcal{P}}(N)$ dividing $q_j$.\\
Then $m>q_j/2$.
\end{Thm}
Note that if for some integer $N$ we could prove that if all continued fraction digits $a_i$ satisfy
$a_i\le 100N$ say and $N|k$, then (\ref{EME}) would be resolved!  
However, for a generic number $\xi\in [0,1]$ that is not a rational, 
one can show that the sequence of $a_i$ is not bounded above.
The Gauss-Kuz'min statistics make this more precise and assert that the probability
that a given term in the continued fraction expansion of a generic $\xi$ is at least $b$, equals $\log_2(1+1/b)$. 
Thus for
a sufficiently large $N$, one expects that $j(N)$ is quite large. This, in combination with the exponential growth
of $q_j$ then ensures a large lower bound for $m$. (The numbers $(\log 2)/2N$ are expected to be generic.)

The conditions b and c are of lesser importance. It seems that condition b is satisfied with probabilty $1/2$.
In practice, sometimes condition a is satisfied, but not b or c, and this leads to a larger lower bound
for $m$.
Condition c is derived using the Moser method, namely by analyzing the equation
\begin{equation}
\label{sta}
\frac{2(3^k-1)(m-1)^k}{2m-3}\equiv -\sum_{\substack{p\mid 2m-3\\p-1\mid k}}\frac1p\pmod1,
\end{equation}
that a solution $(m,k)$ of (\ref{EME}) must satisfy.

We leave it as a challenge to experts in the metric theory of continued fractions to determine the expected
value of $q_{j(N)}$ on replacing $({\log 2})/(2N)$ above by a generic number $\xi$. Gallot et al. expect
that conditions a and b lead to $E(\log q_{J(N)}(\xi))\sim c_1 N$ and, taking into account also condition c,
$E(\log q_{J(N)}(\xi))\sim c_2 N\log^{\beta} N$ for some positive constants $c_1,c_2$ and $\beta$.

Crucial in applying the result is a very good algorithm to determine $\log 2$ with many decimals of
accuracy. Indeed, it is a well-known result of Lochs, that if one knows a generic number up to $n$ decimal 
digits, the one can accurately compute approximately  0.97n continued fraction digits.  For example, knowing
1000 decimal digits of $\pi$ allows one to compute 968 continued fraction digits.

\indent Applying Theorem~\ref{main3} with $N=2^8\cdot 3^5\cdot 5^3$
or $N=2^8\cdot 3^5\cdot 5^4$, and using that
$N|N_2$ and $N_2|k$, Gallot et al. obtained the current world record:
\begin{Thm}
\label{YVES}
If an integer pair $(m,k)$ with $k\ge 2$ satisfies \eqref{EME}, then
$$m > 2.713\,9 \cdot 10^{\,1\,667\,658\,416}>10^{10^9}.$$
\end{Thm}
Gallot et al.~argue that, assuming one can compute $\log 2$ with arbitrary precision, applying Theorem \ref{main3}
with $N=N_2$ should give rise to $m>10^{10^{400}}$.

Interestingly, the results obtained by invoking Bernoulli numbers (`arithmetic') and analysis seem to be completely 
unrelated (`the arithmetic does not feel the analysis'). This strongly suggests that the Erd\H{o}s-Moser conjecture
ought to be true.

\section{A new result}
\label{sec:new}
This section focuses on new research; familiarity with the theory of divisors
of second order sequences is helpful. The reader is referred to Ballot \cite{Ballot} or Moree \cite{Artin} for
more introductory accounts.\\
\indent Let $S$ be an infinite sequence of positive integers. We say that a prime $p$ divides the sequence if it divides
at least one of its terms. Here we will be interested in the sequence $S_2:=\{2^{2e+1}+1\}_{e=0}^{\infty}$. 
It can be shown that $p>2$ divides $S_2$ iff 
ord$_2(p)\equiv 2~({\rm mod~}4)$, with ord$_g(p)$ (with $p\nmid g$) the smallest positive integer $t$ such that
$g^t\equiv 1~({\rm mod~}p)$. 
The set of these primes is known to have natural
density $7/24$ \cite{ak+bk}. Furthermore, if ord$_2(p)\equiv 2~({\rm mod~}4)$ then 
\begin{equation}
\label{bloepie}
p|2^{2e+1}+1{\rm ~iff~}2e\equiv {{\rm ord}_2(p)\over 2}-1 ~({\rm mod~}{\rm ord}_2(p)).
\end{equation}
In some coding theoretical work the sequence $S_2$ and its variants play an important role, as in \cite{DMS, JLX} 
and similarly in the  study of the Stufe of cyclotomic fields \cite{FGS, ak+bk} and the
study of Fermat varieties \cite{Nowak, Shioda}.\\ 
\indent If $m+2$ is coprime with $S_2$, then from (\ref{m+2}) and $2|k$ we can infer
a fifth identity of Moser type, (\ref{newbee}). This then leads to
$m>10^{10^{11}}$ for such $m$. 
We now consider the situation in greater detail.
\begin{Thm}
\label{thm:m+2}
Let $N\equiv 0~({\rm mod~}24)$ be an arbitrary integer. Suppose that $(m,k)$ is a solution of {\rm (\ref{EME})} with 
$$k\ge 2,~N|k{\rm ~and~}m<10^{10^{11}},$$ then
$m+2$ has a prime divisor $p>3$ such that:\\
{\rm 1)} $({\rm ord}_2(p),N)=2$;\\
{\rm 2)} $k\equiv {{\rm ord}_2(p)\over 2}-1 ~({\rm mod~}{\rm ord}_2(p))$.\\
In case $m\equiv 2~({\rm mod~}3)$, we can replace $10^{10^{11}}$ by $10^{10^{16}}$.
In case $N=N_2$ we have $p\ge 2\,099$.
\end{Thm}
\indent We first prove a corollary.
\begin{cor}
\label{corcor}
Suppose every prime divisor $p>3$ of $m+2$ satisfies $p\equiv 5,7~({\rm mod~}8)$. Then
\begin{equation}
\label{gevalletjes}
m\ge
\begin{cases}
10^{10^{16}} & {\rm if~}3\nmid m+2;\\
10^{10^{11}} & {\rm if~}3\mid m+2.
\end{cases}
\end{equation}
\end{cor}
\noindent {\it Proof}. Using the supplementary law of quadratic reciprocity, $({2\over p})=(-1)^{(p^2-1)/8}$, one sees
that if $p\equiv 5,7~({\rm mod~}8)$, then ord$_2(p)\not\equiv 2~({\rm mod~}4)$. Thus condition 1
is not satisfied, as for it to be
satisfied we must have ${\rm ord}_2(p)\equiv 2~({\rm mod~}4)$.\qed\\ 

\noindent Put $$P(N)=\{p>3: ({\rm ord}_2(p),N)=2\}.$$
\noindent We will study the set $P(N)$ in greater detail with the ultimate goal of studying the
$N$-good integers, that is the odd integers $n$ having no prime divisors in $P(N)$. Note that in the proof
of Corollary \ref{corcor}, we established that integers composed only of primes $p\equiv 5,7~({\rm mod~}8)$ are
$N$-good (with $24|N$).
\begin{cor}
\label{corcor2}
Let $N\equiv 0~({\rm mod~}24)$ be an arbitrary integer. If $(m,k)$ satisfies {\rm (\ref{EME})}, $N|k$ and $m+2$
is $N$-good, then $m$ 
satisfies inequality {\rm (\ref{gevalletjes})}.
\end{cor}
\indent If $p$
is to be in $P(N)$, then $p\equiv 1~({\rm mod~}8)$ or $p\equiv 3~({\rm mod~}8)$. 
In the
latter case we have ord$_2(p)\equiv 2~({\rm mod~}4)$. In the former case it is not
necessarily so that ord$_2(p)\equiv 2~({\rm mod~}4)$, and numerically there is a strong preponderance
of primes $p\equiv 3~({\rm mod~}8)$ in $P(N)$. Indeed, we have the following result.
\begin{Lem}
\label{allesisrelatief}
The relative density of primes 
$p\equiv 1~({\rm mod~}8)$ satisfying ord$_2(p)\equiv 2~({\rm mod~}4)$ within the set of
primes $p\equiv 1~({\rm mod~}8)$ is
$1/6$.
\end{Lem}
\noindent {\it Proof}. We have seen that if ord$_2(p)\equiv 2~({\rm mod~}4)$, then  
$p\equiv 1,3~({\rm mod~}8)$. If $p\equiv 3~({\rm mod~}8)$, then ord$_2(p)\equiv 2~({\rm mod~}4)$.
From this, the fact that $\delta({\rm ord}_2(p)\equiv 2~({\rm mod~}4))=7/24$ and the prime number
theorem for primes in arithmetic progression, we infer that the density of primes
$p\equiv 1~({\rm mod~}8)$ such that ord$_2(p)\equiv 2~({\rm mod~}4)$ equals ${7\over 24}-{1\over 4}={1\over
24}$. The sought for relative density is then ${1\over 24}/{1\over 4}={1\over 6}$. \qed\\

\noindent Thus if $p\equiv 3~({\rm mod~}8)$, then ord$_2(p)\equiv 2~({\rm mod~}4)$ and if
$p\equiv 1~({\rm mod~}8)$, then in $1/6$-th of the cases we have
ord$_2(p)\equiv 2~({\rm mod~}4)$.

A further observation concerning the set $P(N)$ is related to Sophie Germain primes.
A prime $q$ such that $2q+1$ is a prime, is called
a {\it Sophie Germain prime}. Let $q_M$ denote the largest prime factor of $M$.
\begin{Lem}
Let $N\equiv 0~({\rm mod~}24)$ be an arbitrary integer.
If $q$ is a Sophie Germain prime, $q\equiv 1~({\rm mod~}4)$ and $q$ and $N$ are coprime, then
$p=2q+1\in P(N)$.
\end{Lem}
\noindent {\it Proof}. The assumptions imply that $({2\over p})=-1$ and since $p>3$ we infer
that ord$_2(p)=2q$. Since $({\rm ord}_2(p),N)=(2q,N)=2$ we are done. \qed\\

\noindent There are $42$ primes $p$ in $P(N_2)$ not exceeding 10~000. Of those $7$ primes
$p$ are such that $(p-1)/2$ is not Sophie Germain, the 
smallest one being $7\,699$. However, the Sophie Germain primes have
natural density zero, whereas as we shall see $P(N)$ has positive natural density.\\
\indent Given a rational number $g$ such that $g\not\in \{-1,0,1\}$, the natural density
$\delta_g(d)$ of the set of primes $p$ such that the order of $g~({\rm mod~}p)$ is
divisible by $d$ is known to exist and can be computed; see e.g. Moree \cite{polen}.
Using inclusion and exclusion one then finds that the set $P(N)$ has natural density
$$\delta(N)=\sum_{d|N_0}(\delta_2(2d)-\delta_2(4d))\mu(d),$$
where $N_0$ is the product of the odd prime divisors dividing $N$ and
$\mu$ denotes the M\"obius function. By Moree \cite[Theorem 2]{polen} we
then find that, for odd $d$,
$$\delta_2(2d)-\delta_2(4d)={7\over 24}\prod_{p|d}{p\over p^2-1},$$
and hence
$$\delta(N)={7\over 24}\sum_{d|N_0}\mu(d)\prod_{p|d}{p\over p^2-1}={7\over 24}\prod_{p|N_0}\left(1-{p\over p^2-1}\right),$$
where we used that a multiplicative function $f$ satisfies 
$$\sum_{d|N_0}\mu(d)f(d)=\prod_{p|N_0}(1-f(p)).$$
Taking $N=N_2$ one finds that
$$\delta(N_2)={7\over 24}\prod_{2<p\le 1000}\Big(1-{p\over p^2-1}\Big)\approx 0.043\,578\,833\cdots$$
By a result of Wiertelak, quoted as Theorem 1 in Moree \cite{polen}, we have
$$\sum_{p\le x,~p\not\in P(N)}1=(1-\delta(N)){x\over \log x}+O_N({x\over \log^2 x}),$$
where the implicit constant may depend on $N$. From this result and \cite[Proposition 4]{Martin}, we then
infer that asymptotically the number of integers $n\le x$ that
are $N$-good, $N_G(x)$, satisfies 
$$N_G(x)\sim c_Nx \log^{-\delta(N)}x,$$ where 
$$c_N={1\over \Gamma(1-\delta(N))}\lim_{x\rightarrow \infty}\prod_{p\le x}\left(1-{1\over p}\right)^{1-\delta(N)}\left(1-{\chi_N(p)\over p}\right)^{-1},$$
with $\chi_N(p)=0$ if $p=2$ or $p$ is in $P(N)$ and 1 otherwise. 
(As usual $\Gamma$ denotes the Gamma-function.)
Taking $N=N_2$, a computer calculation
suggests that $c_{N_2}\approx 0.54$.\\

Now if we have a sequence of random integers $n_j$ growing roughly as $e^{\beta j}$ for some constant $\beta>0$, the integer
$n_j$ is $N$-good with probability $c_N\log^{-\delta(N)}n_j\approx c_N(\beta j)^{-\delta(N)}$. The expected
number of $N$-good $n_j$ with $j\le x$ is then approximately
$$c_N\sum_{j\le x} (\beta j)^{-\delta(N)}\sim c_N{(\beta x)^{1-\delta(N)}\over (1-\delta(N))\beta}.$$
The result that $2k/(2m-3)$ is a convergent $p_j/q_j$ of $\log 2$ with $j$ even and the result of L\'evy \cite{Paul}
that for a generic $\xi\in [0,1]$ that is not a rational
$$\lim_{j\rightarrow \infty}{\log q_j(\xi)\over j}={\pi^2\over 12\log 2}\approx 1.18,$$
leads us to expect that the sequence $m_j$ of potential solutions $(k_j,m_j)$ to (\ref{EME}) coming from this result, 
is of exponential growth. Thus of the potential solutions $(m_j,k_j)$ with
$j\le x$, one expects about $x^{1-\delta(N_2)}$, that
is roughly $x^{0.96}$, to be $N_2$-good. 
For those (\ref{gevalletjes}) holds with $m=m_j$.
Thus if there would be say $10^{10}$ potential solutions
with $m\le 10^{10^{11}}$, then one expects roughly $3\cdot 10^9$ to be $N_2$-good, and those can
be excluded by Corollary \ref{corcor2}.\\

\noindent {\tt Remark}. Given positive integers $a,b,c,d$, the density of primes $p\equiv c~({\rm mod~}d)$ such that
$p|\{a^e+b^e\}_{e=0}^{\infty}$ is known; see Moree and Sury \cite{MSury}. Since $S_2=\{2\cdot
4^e+1\}_{e=0}^{\infty}$, that result cannot be applied to establish Lemma \ref{allesisrelatief}.\\

\noindent {\it Proof of Theorem} \ref{thm:m+2}. The idea of the proof is to show that if for every 
prime divisor $p>3$ of $m+2$ at
least one the conditions 1 or 2 is not satisfied, then the identity   
\begin{equation}
\label{newbee}
\sum_{p|m+2}{1\over p}+{3\over m+2}\equiv 0~({\rm mod~}1)
\end{equation}
holds. Using this we then show that $m$ is bigger than the bound in the theorem; this is a  contradiction. As
usual we make heavy use of the fact that $k$ must be even.\\
\indent We start with the equation
\begin{equation}
\label{m+2}
2^{k+1}+1\equiv -\sum_{p-1|k,~p|m+2}{(m+2)\over p}~({\rm mod~}m+2),
\end{equation}
found on noting that $S_k(m+2)=2m^k+(m+1)^k\equiv 2^{k+1}+1~({\rm mod~}m+2)$ and on invoking
Theorem \ref{thm:clausen}.
Suppose that $p|m+2$. The idea
is to reduce (\ref{m+2}) modulo $p$ (except if $p=3$, then we reduce modulo $9$).\\
\indent If $p=3$, then using $6|k$ we see that $2^{k+1}+1\equiv 3~({\rm mod~}9)$, and
we infer that $3^2||m+2$, that is we must have $m\equiv 7,16~({\rm mod~}27)$. Next assume $p>3$.\\
\indent First assume that ord$_2(p)\not\equiv 2~({\rm mod~}4)$. 
Then $p$ does not divide $S_2$. Thus the right hand side
of (\ref{m+2}) is non-zero modulo $p$, and this implies that $p-1|k$ and $p^2\nmid m+2$ and
hence $2^{k+1}+1\equiv 3~({\rm mod~}p)$.\\
\indent Next assume that ord$_2(p)\equiv 2~({\rm mod~}4)$, and condition 1 is not 
satisfied. Then ord$_2(p)$ and $N$ have an odd prime factor in common, and by
(\ref{bloepie}) (with $e=k/2$) we get a contradiction to the assumption $N|k$.\\
\indent Finally, assuming that condition 1 is satisfied but not condition 2, the right hand side
of (\ref{m+2}) is non-zero modulo $p$, and the same conclusion as before holds.
By the Chinese remainder theorem we then infer that $2^{k+1}+1\equiv 3~({\rm mod~}m+2)$, and
hence from (\ref{m+2}) we see that (\ref{newbee}) holds.\\
\indent Put $M_3=(m^2-1)(4m^2-1)(m+2)$. By part 2 of Theorem \ref{main} we infer that amongst 
the numbers $m-1,m+1,m+2,2m-1,2m+1$, no prime $p\ge 7$ occurs more than once as
divisor, the prime 2 occurs precisely twice, the prime 3 at most 3 times and the prime 5 at most two times.
Using this, we obtain on adding Moser's equations (\ref{M8}), (\ref{M13}), (\ref{M16}) and (\ref{M18})
to (\ref{newbee}):
\begin{equation}
\label{improvement}
\sum_{p|M_3}{1\over p}+{1\over m-1}+{2\over m+1}+{2\over 2m-1}+{4\over 2m+1}+{3\over m+2}\ge {109\over 30},
\end{equation}
where $${109\over 30}=5-{1\over 2}-{2\over 3}-{1\over 5}=3.633\,333\,333\,333\cdots $$ Using the estimate
$$\sum_{p\le x}{1\over p}<\log \log x + 0.2615+{1\over \log^2 x}{\rm ~for~}x>1,$$
due to Rosser and Schoenfeld \cite[(3.20)]{RS1}, we find that
$\sum_{p\le \beta}1/p<3.633\,32$ with $\beta=4.33\cdot 10^{12}$. {}From another paper
by the same authors \cite{RS2} we have
$$|\theta(x)-x|<{x\over 40\log x},~x\ge 678\,407.$$
Hence
$$\log(4m^5)>\log(N_3)>\log\prod_{p\le \beta}p=\theta(\beta)>.999\beta ,$$
from which we infer that $m\ge 10^{10^{11}}$.\\
\indent In case $m\equiv 2~({\rm mod~}3)$ there are precisely two of the five terms
$m-1,m+1,2m-1,2m+1$ and $m+2$ divisible by $3$, and in (\ref{improvement}) we can replace
$109/30$ by $109/30+1/3=119/30=3.966\,666\cdots$. In that case we can take $\beta=4.425\cdot 10^{17}$ and this leads
to $m\ge 10^{10^{16}}$.\\
\indent The smallest two primes in $P(N_2)$
are $2\,027$ and $2\,099$. For $p=2\,027$ we can actually show that condition 2 is not satisfied. To
this end we must show that $k\not\equiv 1\,012~({\rm mod~}2\,026)$. Computation shows that
$(1\,012,6\,079),~(3\,038,6\,079)$ and $(5\,064,6\,079)$ are helpful pairs. 
By Lemma \ref{helpful} it then follows that $k\not\equiv 1\,012~({\rm mod~}2\,026)$.
The smallest prime that possibly
satisfies both condition 1 and 2 is hence $2\,099$. \qed\\

\noindent{\tt Remark 1}. We leave it as an exercise to the reader to show that (\ref{gevalletjes}) can be refined to
\begin{equation}
\label{gevalletjes2}
m\ge
\begin{cases}
10^{10^{20}} & {\rm if~}3\nmid m+2;~5\nmid m+2\\
10^{10^{16}} & {\rm if~}3\nmid m+2;~5\mid m+2\\
10^{10^{14}} & {\rm if~}3\mid m+2;~5\nmid m+2\\
10^{10^{11}} & {\rm if~}3\mid m+2;~5\mid m+2. 
\end{cases}
\end{equation}
In the same vein one can show that if $(m,k)$ satisfies (\ref{EME}), $k\ge 2$ and
$m\equiv \pm 1~({\rm mod~}15)$ or $m\equiv \pm 1~({\rm mod~}21)$, then
$m\ge 10^{10^{20}}$. If, e.g., $m\equiv 1~({\rm mod~}15)$, then the sum in the left hand side of (\ref{M8}) exceeds
1, so must be at least two. We infer that (\ref{M19}) holds with $3.1666...$ replaced by $4.1666....$. This then
leads to $m\ge 10^{10^{20}}$. The remaining cases are similar (they all lead to (\ref{M19}) with $3.1666...$ 
replaced by $4.1666....$).\\

\noindent{\tt Remark 2}. Using the methods from Bach et al. \cite{BKS}, it
should be possible to compute the largest $\beta$ such that $\sum_{p\le \beta}1/p<109/30$, respectively $119/30$ exactly. They found, e.g., that the 
prime $p_0=180\,124\,123\,005\,660\,046\,7$ is the largest one such that $\sum_{p\le p_0}1/p<4$.

\section{The generalized Erd\H{o}s-Moser conjecture}
\label{eight}
\noindent The Erd\H{o}s-Moser conjecture has the following generalization.
\begin{Con}
There are no integer solutions $(m,k,a)$ of
\begin{equation}
\label{EMEg}
1^k+2^k+\ldots+(m-1)^k=am^k
\end{equation}
with $k\ge 2$, $m\ge 2$ and $a\ge 1$.
\end{Con}
\noindent In this direction the author proved in 1996 \cite{Oz} that (\ref{EMEg}) has no integer
solutions $(a,m,k)$ with $k>1$ and $m<\max(10^{10^6},a\cdot 10^{22})$. With
the hindsight of more than 10 years this can be improved.
\begin{Thm}
\label{provo}
The equation {\rm (\ref{EMEg})} has no integer solutions $(a,m,k)$ with
$$k\ge 2,~m<\max \big(10^{9\cdot 10^6},a\cdot 10^{28}\Big).$$
\end{Thm} 
\noindent {\it Proof}. 
(In this proof references to propositions and lemmas are exclusively to those
in \cite{Oz}.)
The Moser method yields that $2|k$ and gives the following four inequalities
\begin{equation}
\label{P11}
\sum_{p-1|k\atop p|m-1}{1\over p}+{a\over m-1}\ge 1,~~
\sum_{p-1|k\atop p|m+1}{1\over p}+{a+1\over m+1}\ge 1.
\end{equation}
\begin{equation}
\label{P12}
\sum_{p-1|k\atop p|2m-1}{1\over p}+{2a\over 2m-1}\ge 1,~~
\sum_{p-1|k\atop p|2m+1}{1\over p}+{2(a+1)\over 2m+1}\ge 1.
\end{equation}
Since $p|m$ implies $p-1\nmid k$ (Proposition 9), we infer that
$(6,m)=1$. Using this we see that $M_1=(m^2-1)(4m^2-1)/12$ is an even integer.
Since no prime $>3$ can divide more than one of the numbers $m-1,m+1,2m-1$ and $2m+1$, and
since 2 and 3 divide two of these numbers, we find on adding the inequalities that
$$\sum_{p-1|k,~p|M_1}{1\over p}+{a\over m-1}+{a+1\over m+1}+{2a\over 2m-1}+{2(a+1)\over 2m+1}
\ge 4-{1\over 2}-{1\over 3}=3{1\over 6}.$$
Using that $a(k+1)<m<(a+1)(k+1)$ (Proposition 2), we see that in the latter equation
the four terms involving $a$ are bounded above by $6/(k+1)$. Since $k\ge 10^{22}$ 
(Lemma 2), we can proceed as in the proof of Theorem \ref{main} and find the
same bound for $m$, namely $m > 1.485\cdot 10^{9321155}$.\\
\indent Earlier it was shown that if $k>1$, then $k\ge 10^{22}$. To this end
Proposition 6 with $C=3.16$, $s=664\,579=\pi(10^7)$ and $n$ the 200-th highly composite
number was applied. Instead we apply it with $C=3{1\over 6}-10^{-10}$,
$s=4\,990\,906$ and $n$ the 259-th composite number $c_{250}$ (this has the 
property that the number of divisors of $c_{259}<s$, whereas the number of divisors
of $c_{260}$ exceeds $s$). Since
$n=c_{259}>5.583\,4\cdot 10^{27}$ it follows that $k\ge 2n>10^{28}$. Since $m>a(k+1)$, the
proof is completed. \qed\\ 

\noindent {\tt Remark 1}. The above proof shows that if (\ref{EMEg}) has a solution
with $k\ge 2$, $m\ge 2$ and $a\ge 1$, then $m$ must be odd. An easy reproof of this
was given by MacMillan and Sondow \cite{MS2}.\\

\noindent {\tt Remark 2}. The reader might wonder whether the method of Gallot et al. can be applied here as well to
break the $10^{10^7}$ barrier. For a fixed integer $a$ this is possible if one manages to establish that 
$N|k$ with $N$ large enough. Gallot et al. showed that $2k/(2m-2a-1)$ is a convergent with even index
of $\log(1+1/a)$ for $m$ large enough. For a given $a$ this can be made effective. Establishing that $N|k$ along the lines
of Section \ref{sec:MRU} is not always possible (see the last paragraph of that section).\\

\noindent {\tt Challenge}: Reach the benchmark $10^{10^7}$ in Theorem \ref{provo}.

\section{The Kellner-Erd\H{o}s-Moser conjecture}
\label{sec:kell}
Kellner \cite{K} conjectured that if $k,m$ are positive integers with $m\ge 3$, the
ratio $S_k(m+1)/S_k(m)$ is an integer iff $(k,m)\in \{(1,3),(3,3)\}$. Noting that
$S_k(m+1)=S_k(m)+m^k$ one easily observes that this conjecture is equivalent with the
following one.
\begin{Con}
\label{kellie}
We have $aS_k(m)=m^k$ iff $(a,k,m)\in \{(1,1,3),(3,3,3)\}$.
\end{Con}
\noindent If this conjecture holds true, then obviously so does the Erd\H{o}s-Moser conjecture.\\
\indent It is easy to deal with the case $m=3$.  Then we must have $a(1+2^k)=3^k$, and hence
$a=3^e$ for some $e\le k$. It follows that $1+2^k=3^{k-e}$. This Diophantine equation was
already solved by the famous medieval astronomer Levi ben Gerson (1288-1344), alias Leo Hebraeus, who
showed that 8 and 9 are the only consecutive integers in the sequence of powers of 2 and 3, see
Ribenboim \cite[pp. 124-125]{Ribenboim}. This leads to the solutions $(e,k)\in \{(0,1),(3,1)\}$ and hence
$(a,k,m)\in \{(1,1,3),(3,3,3)\}$. Next assume that $m\ge 4$ and $k$ is odd. Then by Theorem
\ref{thm:clausen} we find that $m(m-1)/2$ divides $m^k$, which is impossible. We infer that to 
establish Conjecture \ref{kellie}, it is enough to establish Conjecture \ref{kellie2}, where
$${\mathcal A}=\{a\ge 1: aS_k(m)=m^k{\rm ~has~a~solution~with~}2|k,~k\ge 2,~m\ge 4\}.$$
\begin{Con}
\label{kellie2}
The set $\mathcal A$ is empty.
\end{Con}
The next result shows that if $a\equiv 2~({\rm mod~}4)$ or $a\equiv 3,6~({\rm mod~}9)$, then
$a\not \in \mathcal A$.
\begin{Thm} Let $k\ge 2$ be even. Suppose that $q|a$ is a prime such that $q^2\nmid a$ and
$q-1|k$. Then $aS_k(m)\ne m^k$ and hence $a\not\in \mathcal A$.
\end{Thm}
{\it Proof}. Suppose that $aS_k(m)=m^k$. Let $q^e||m$. Note that $e\ge 1$. 
Using Theorem \ref{thm:clausen} we find that 
$S_k(m)\equiv {m\over q^e}S_k(q^e)\equiv -{m\over q}~({\rm mod~}q^e)$. Now we consider the
identity $aS_k(m)=m^k$ modulo $q^{e+1}$ and find
$-a{m\over q}\equiv m^k\equiv 0~({\rm mod~}q^{e+1})$, contradicting $q^{e+1}||am$. It
follows that $aS_k(m)\ne m^k$.\qed\\

Note that if $a\not\in \mathcal A$, then the equation $aS_k(m)=m^k$ can be solved completely. The
author is not aware of earlier `naturally' occurring Erd\H{o}s-Moser type equations that can be
solved completely. He expects that further values of $a$ can be excluded and might come back to this in a 
future publication. \\

\noindent {\tt Acknowledgement}. 
Part of this article was written whilst I had 4 interns (Valentin Buciumas, Raluca Havarneanu, Necla Kayaalp and Muriel Lang) studying variants
of the Erd\H{o}s-Moser equation. Raluca and Valentin found a (reparable) mistake in
Moser's paper and Muriel showed that $2027$ is the smallest prime in $P(N_2)$, but does 
not satisfy condition 
2 of Theorem \ref{thm:m+2}. I thank them all for their
questions, comments and cheerful presence. Paul Tegelaar provided some helpful comments on an earlier version. Jonathan Sondow I thank for helpful e-mail correspondence. 
This note profited a lot from corrections by Julie Rowlett (a native English speaker).
Particular thanks are due to the referee for 
many 
very detailed and constructive comments.\\
\indent The academic year 1994/1995 the author spent as a postdoc of Alf van der Poorten at
Macquarie. Alf told me various times it would be so nice if mathematicians could be less serious
in their mathematical presentation, e.g. talk about a `troublesome double sum', if there are
smooth numbers, then also consider hairy numbers, etc.. In this spirit, I `spiced up' my
initial submission of \cite{Oz}, the red pencil of the referee was harsh though, but somehow
the words `mathemagics' and `rabbits' survived. Rabbits are difficult to suppress, and least
of all mathemagical ones. So I am happy they are back full force in the title of \cite{rabbits}. 
Also whilst at Macquarie, thanks to questions by then visitor Patrick Sol\'e, I got into the
study of divisors of $a^k+b^k$, not realizing there is a connection with the Erd\H{o}s-Moser equation (as
the present article shows).\\
\indent Had Alf learned that the present record for solutions for EM is based on a continued 
fraction expansion (of $\log 2$), I am sure he would have been pleased.\\
\indent I had a wonderful year in Australia and will be always grateful to Alf for having made
that possible.

\end{document}